\documentstyle{article}

\newcommand{\R}{{\bf R}}
\newcommand{\T}{{\bf T}}

\newcommand{\Q}{{\bf Q}}
\newcommand{\N}{{\bf N}}
\newcommand{\Z}{{\bf Z}}

\newcommand{\CP}{{\bf CP}}

\newcommand{\eps}{{\varepsilon}}

\newcommand{\length} {{\rm length}}
\newcommand{\MS}{{\medskip}}

\newcommand{\NI}{{\noindent}}

\newcommand{\QED}{\hfill$\Box$\medskip}

\begin{document}

\title{Hamiltonian loops from the ergodic point of view}
\author{
Leonid
Polterovich\thanks{Supported by the United States - Israel
Binational Science Foundation grant 94-00302
} \\ School of Mathematical Sciences\\
Tel Aviv University \\ 69978 Tel Aviv, Israel\\ (polterov@math.tau.ac.il)}

\date{} 
\maketitle

\NI

Let $G$ be the group of Hamiltonian diffeomorphisms of a closed
symplectic manifold $Y$.
A loop $h:S^1 \to G$ is called strictly ergodic if 
for some irrational number $\alpha$
the associated skew product map $T:S^1\times Y \to S^1 \times Y$
defined by $T(t,y) = (t + \alpha, h(t)y)$ is strictly ergodic.
In the present paper we address the following question. Which
elements of the fundamental group of $G$ can be represented
by strictly ergodic loops? We prove existence of contractible
strictly ergodic loops for a wide class of symplectic manifolds
(for instance for simply connected ones). Further,
we find a restriction
on the homotopy classes of smooth strictly ergodic loops in the framework of
Hofer's bi-invariant geometry on $G$. Namely, we prove that
their asymptotic Hofer's norm must vanish. This result provides
a link between ergodic theory and symplectic 
topology. 

\NI

\bigskip

{\bf 1. Introduction and results}

\bigskip

{\it 1.1 Hamiltonian loops as dynamical objects}

\medskip

Let $(Y,\Omega)$ be a closed symplectic manifold, and
let $G$ be its group of Hamiltonian diffeomorphisms.
\footnote{
Recall that $G$ consists of all symplectomorphisms of $(Y,\Omega)$
which can be included into a time-dependent
Hamiltonian flow on $Y$.
}
Given an irrational number $\alpha$ and a smooth loop
$h: S^1 \to G$, one can define {\it a
skew product map} $T_{h,\alpha}:S^1\times Y \to S^1 \times Y$
by $T_{h,\alpha}(t,y) = (t + \alpha, h(t)y)$.
The purpose of the present paper is to relate geometry and topology
of Hamiltonian loops with dynamics of associated skew products.

The dynamical property we consider is the strict ergodicity.
Recall that a homeomorphism $T$ of a compact topological space
$X$
is called {\it strictly ergodic} if it has precisely one invariant
Borel probability measure, say $\mu$ , and moreover this measure is positive
on non-empty open subsets. Strictly ergodic homeomorphisms are ergodic,
and have a number
of additional remarkable features.
We mention one of them which plays
a crucial role below. Namely, given such a  $T$ and an arbitrary
continuous function $F$ on $X$, the time averages
${1\over N}\Sigma_{i=0}^{N-1}F(T^ix)$ converge uniformly to the space average
$\int_X Fd\mu$, and in particular converge
{\it for all} $x \in X$. Note that for general ergodic
transformations such a convergence takes place only for {\it almost
all} $x$. The contrast between "all" and "almost all" becomes especially
transparent when one notices that there are pure topological obstructions
to the strict ergodicity. For instance, the 2-sphere admits no strictly
ergodic homeomorphisms. Indeed the Lefschetz theorem implies that
every homeomorphism of $S^2$ has either a fixed point, or a periodic
orbit of period $2$ and we see that the invariant measure
which is concentrated on such an orbit contradicts the
definition of strict ergodicity. In 1.3 below we  describe a more
sophisticated obstruction to the strict ergodicity which comes from
symplectic topology.

We say that 
a loop $h:S^1 \to G$ is {\it strictly
ergodic} if for some $\alpha$
the corresponding skew product map $T_{h,\alpha}$ is strictly ergodic.
With this language our central question can be formulated as follows.

\medskip

\NI
{\bf Question 1.1.A.} Which homotopy classes $S^1 \to G$ can be
represented by strictly ergodic loops? 

\medskip

Here is an example where one gets a complete answer to this question.
Let $Y$ be the blow up of the complex
projective plane $\CP^2$ at one point. Choose a K\"ahler symplectic
structure $\Omega$
on $Y$ which integrates to 1 over a general line and to $1\over 3$
over the exceptional divisor. The periods of the
symplectic form are chosen in such a way that its
cohomology class is a multiple of the first Chern class of $Y$.
One can easily see that $(Y,\Omega)$ admits an
effective Hamiltonian action
of the unitary group $U(2)$, in other words there exists a
monomorphism $i: U(2) \to G$. The fundamental group of $U(2)$
equals $\Z$. 
It was
proved recently by Abreu and McDuff [AM] that 
the inclusion $\pi_1(U(2)) \to \pi_1(G)$ is an isomorphism,
and thus $\pi_1(G) = \Z$. As far as we know this is the simplest
example of a symplectic manifold with $\pi_1(G)=\Z$.

\proclaim Theorem 1.1.B.
In this situation, the trivial class $0 \in \Z$ is the only  
one which can be represented by a strictly ergodic loop.

The proof (see 1.3 below) is based on two general results
on existence and non-existence of strictly ergodic loops.

\medskip

{\it 1.2 An existence result}

\medskip

Assume in addition that
the group $G$ of all Hamiltonian diffeomorphisms of a closed
symplectic manifold  $(Y,\Omega)$
is $C^{\infty}$-closed in ${\rm Diff}(Y)$.

\NI

\proclaim Theorem 1.2.A. Under this assumption there
exists a contractible smooth strictly ergodic Hamiltonian loop.

\medskip

\NI
{\bf Remark. } 
The assumption above holds
for a wide class of symplectic manifolds, for instance when $H^1(Y,\R) = 0$,
or when the cohomology class of the symplectic form is rational.
(And thus it holds for the blow up of $\CP^2$ considered in 1.1 above).
However
it is still unclear whether it is valid for all closed symplectic manifolds.
This long standing problem is known as the
Flux conjecture, and we refer the reader to
[LMP] for more discussion and recent results.
On the other hand it sounds likely that in our situation
this assumption plays a technical role only and can be removed,
but I have not checked the details (see the footnote
in 5.3 below; I am thankful to D. McDuff
for illuminating discussions on this issue).

\medskip

Constructions of ergodic and strictly ergodic skew products
associated to a loop in a group acting on a topological space
have a long history (see [AK], [FH], [GW], [He], [N]).
All these constructions are based
on a beautiful and quite counter-intuitive idea to look
for such skew products in the closure of ones with
absolutely trivial dynamical behaviour.
\footnote
{
A different approach based on KAM theory
was used in [E], see also 1.7.B below.
} 
Here is the precise statement in our setting.

Consider the set $C^{\infty}(S^1,G)$  of all smooth loops
$S^1 \to G$ as a subset of the space $C^{\infty}(S^1 \times Y,Y)$,
and endow it with the topology induced by the $C^{\infty}$-topology.
Consider the subset
${\cal R} \subset S^1 \times C^{\infty}(S^1,G)$ which is defined
as the closure of the following subset:
$$\{(\alpha,h)|\; h(t)=g(t+\alpha)^{-1}g(t), \;
{\rm where}\; g \; {\rm runs} \;{\rm over} \; C^{\infty}(S^1,G)
\}.$$

\proclaim Theorem 1.2.B.
The pairs $(\alpha,h)$ such that the transformation
$(t,y) \to (t+\alpha,h(t)y)$
is strictly ergodic form a residual
subset in $\cal R$.

Though in the literature there are plenty of similar statements,
those of them which I found deal either with other groups $G$,
or with usual ergodicity. 
We outline the proof of 1.2.B in \S 2
below, and give full details in \S 3 - \S 5.

Note that the loops $h$ which appear in the definition of
the set $\cal R$ are limits of contractible loops on $G$, and therefore
are contractible. Thus 1.2.B implies 1.2.A.

\medskip

{\it 1.3 An obstruction via Hofer's geometry}

\medskip

Let $Y$ be a closed connected symplectic manifold and let
$G$ be the group of all Hamiltonian diffeomorphisms. 
In 1990 Hofer [H] discovered a remarkable
bi-invariant geometry on $G$.
The development of this geometry has lead to 
a new geometric intuition in dynamical systems (see discussion in [P2]), and
our approach to Question 1.1.A may be considered as one more step
in this direction. Let us introduce the notion of
the (asymptotic) length spectrum in Hofer's geometry which is
relevant for our study of loops of Hamiltonian diffeomorphisms
(see [P1]).

Every 
smooth loop $h: S^1 \to G$ is generated by the unique Hamiltonian
function $H: S^1 \times Y \to \R$ which is normalized as follows:
the integral of $H(t,.)$ over $Y$ vanishes for all $t$.
Define the length of the loop $h$ by
$$\length (h) = \int_0^1 \max_{y \in Y} |H(t,y)| dt.$$
Note that $G$ can be considered as an infinite-dimensional Lie
group whose Lie algebra coincides with the space  $C^{\infty}_0(Y)$
of smooth functions with zero mean. The $L_{\infty}$-norm
on the Lie algebra is invariant under the adjoint action of $G$,
so it defines a bi-invariant Finsler metric on $G$. With this
language the definition of the length above is just the
usual definition of Finsler length. 

Take now an element $\gamma$ of the fundamental group
$\pi_1(G)$. 
Set
$$||\gamma|| = \inf \length (h),$$
where the infimum is taken over all loops $h:S^1 \to G$ representing
$\gamma$. Finally, define {\it asymptotic Hofer's norm} of $\gamma$
as
$$||\gamma||_{\infty} = \lim_{k \to \infty} {1\over k}||\gamma^k||.$$
(The limit exists since the sequence $\{||\gamma^k||\}$ is subadditive).

Using methods of "hard" symplectic topology, one can show that in
some interesting situations
this quantity is non-trivial (see [P1]).

\proclaim Theorem 1.3.A. Let $\gamma \in \pi_1(G)$
be a class represented by a smooth strictly ergodic loop.
Then asymptotic Hofer's norm $||\gamma||_{\infty}$ vanishes.

\medskip

{\bf Proof of Theorem 1.1.B:} Let $Y$ be the monotone blow up of 
$\CP^2$ at one point as in 1.1.B. It follows from 1.2.A
that there exists a contractible strictly ergodic loop. On the
other hand,
it was shown in [P1] that in this case asymptotic Hofer's norm
of every non-trivial element of $\pi_1(G)$ is strictly positive.
Thus there are no non-contractible strictly ergodic loops in view of
Theorem 1.3.A. This completes the proof.
\QED

\NI

The proof of 1.3.A is very simple and we present it immediately
in 1.4.

\medskip

{\it 1.4 Asymptotic shortening of strictly ergodic Hamiltonian loops }

\medskip

Recall that if $h_1(t)$ and $h_2(t)$ are paths of Hamiltonian diffeomorphisms
generated by normalized Hamiltonians $H_1$ and $H_2$ respectively then
the composition $h_2(t) \circ h_1(t)$ is generated by the normalized
Hamiltonian $H_2(t,y) + H_1(t, h_2(t)^{-1}y)$.
Let $h:S^1 \to G$ be a smooth loop of Hamiltonian diffeomorphisms which
defines a strictly
ergodic skew product $T(t,y) = (t+\alpha,h(t)y)$. Let $\gamma$ be the
corresponding element in $\pi_1(G)$. Denote by $H(t,x)$ the normalized
Hamiltonian function generating the loop $h(t)^{-1}$.
Set $h_k(t) = h(t +k\alpha)^{-1}$ and set
$$f_N(t) = h_0(t) \circ ... \circ h_{N-1}(t).$$
In view of the discussion above the loop $f_N$ is generated by the
normalized Hamiltonian function
$$F_N(t,y) = H(t,y)+H(t+\alpha,h_0(t)^{-1}y)+...$$ $$+H(t+(N-1)\alpha,
h_{N-2}(t)^{-1}\circ ... \circ h_0(t)^{-1}y).$$
This expression can be rewritten as follows:
$$F_N(t,y) = \Sigma_{k=0}^{N-1} H\circ T^{k} (t,y).$$
Since $T$ is strictly ergodic and the function $F_N$
has zero mean we conclude that
$${1\over N} \int_0^1 \max_{y \in Y} |F_N(t,y)| dt \to 0,$$
when $N \to \infty$. But the expression on the left hand side is exactly
${1\over N}\length (f_N(t))$. Note now that the loop $f_N(t)$ represents
the element $\gamma^{-N}$. Since $||\gamma^{N}|| = ||\gamma^{-N}||$
we get that
${1\over N}||\gamma^{N}||$
tends to zero when $N \to \infty$. This proves that asymptotic
Hofer's norm of $\gamma$ vanishes.
\QED

\medskip

{\it 1.5 A generalization to sequential systems}

\medskip

We present here a generalization of Theorem 1.3.A which deals with
ergodic properties of so called sequential dynamical
systems. Let $X$ be a compact topological measure space, and let 
$\{T_i\} =(T_1,T_2,...,T_i,...)$ be a sequence of measure-preserving
homeomorphisms.
Such a sequence defines an evolution with discrete
time on $X$. Namely a position of a point $x \in X$ at the time moment
$n \in \N$ is $T^{(n)}(x)$, where here and below we write $T^{(n)}$
for the composition
$T_n\circ...\circ T_1 $. Ergodic properties
of such systems were studied in the literature (see for instance [BB],
as well as an extensive discussion on  random ergodic theorems in [Kr]).
However, we have not found any reference to the next definition
which sounds to us pretty natural. A sequence $\{T_i\}$
is called {\it strictly ergodic} if for every continuous function
$H$ on $X$ with zero mean the time averages
$${1\over N}\Sigma_{i=0}^{N-1} H \circ T^{(i)} $$
converge uniformly to zero. Our convention is that $T^{(0)}$ is 
the identity map.

Let as before $G$ be the group of Hamiltonian diffeomorphisms
of a closed symplectic manifold $(Y,\Omega)$, and $X=S^1\times Y$.
Let $\{\alpha_i\}, \; i \geq 1$ be
an arbitrary sequence of numbers, and let $\{g_i\}$ be an arbitrary
sequence of Hamiltonian diffeomorphisms. Take a smooth loop $h:S^1 \to G$
and consider a sequence $\{T_i\}$ of skew products of the form
$T_i(t,y) = (t+\alpha_i,g_ih(t)y)$. Denote by $\gamma \in \pi_1(G)$ the
element represented by $h$.
In this setting one can generalize
Theorem 1.3.A as follows.

\NI 
{\it If the sequence $T_i$ is strictly ergodic then asymptotic Hofer's
norm of} $\gamma$ {\it vanishes}.

\NI
This can be proved by the following modification of the shortening procedure
described in 1.4 above. Denote by $H(t,y)$ the normalized Hamiltonian
function of $h(t)^{-1}$.

Let $ \phi_0 = {\rm id},
\phi_1, \phi_2,...$ 
be a sequence of transformations from $G$ such that $\phi_i^{-1}\phi_{i-1}=
g_i$ for all $i \geq 1$. Set 
$$h_k(t) = \phi_k h(t+\alpha_0+...+\alpha_k)^{-1}\phi_k^{-1},$$
where $\alpha_0 = 0$.
Consider a new loop $$f_N(t) = h_0(t)\circ...\circ h_{N-1}(t).$$
It is easy to see that $f_N$ is generated by the normalized
Hamiltonian
$$F_N = \Sigma_{i=0}^{N-1} H \circ T^{(i)},$$
and this loop represents the class $\gamma^{-N}$. Now exactly the same
argument as in 1.4 completes the proof.

\medskip

{\it 1.6 An application to Hofer's geometry}

\medskip

I do not know {\it the precise} value of $||\gamma||_{\infty}$ in any
example where this quantity is strictly positive (for instance, 
for the blow up of $\CP^2$ in 1.1 above). The difficulty
is as follows. In all known examples where 
Hofer's norm $||\gamma||$ can be computed precisely
there exists a closed  loop $h(t)$ which minimizes the length
in its homotopy class (that is {\it a minimal closed geodesic}).
It turns out however that every minimal closed geodesic
looses minimality after a suitable number of iterations.
In other words the loop $h(Nt)$ can be shortened provided $N$ is
large enough. The proof of this statement is based on a shortening
procedure described in the previous section and goes as follows.
In the notations of 1.5, take $\alpha_i=0$ for all $i$. 
Set $a(t) =
\max_{y \in Y} |F_N(t,y)|$ and $b(t) = N \max_{y \in Y}|H(t,y)|$.
It is very easy to
choose
a sequence $g_1,....,g_N$ in such a way that  $a(0) < b(0)$.
Since $a(t) \leq b(t)$ for all $t$, we get that
$\int_0^1 a(t) dt < \int_0^1 b(t) dt$,
and this completes the argument.
We conclude that {\it if a class $\gamma \in \pi_1(G)$ is represented
by a minimal geodesic then $||\gamma||_{\infty}$ is strictly less than
$||\gamma||$.} 

Let us complete this section with few remarks on curve shortening
in Hofer's geometry.
The first shortening procedure is due to Sikorav [Si].
Further progress was made by Ustilovsky [U],
Lalonde - McDuff [LM], and in a joint paper
with Bialy [BP]. Our procedures in 1.4 and 1.5 are closely related
to these developments. In particular, in [BP] we asked a question about
the role of Birkhoff's ergodic sums in Hofer's geometry. The results above
can be considered as a sort of answer.

\vfill\eject

{\it 1.7 Remarks and open problems}

\medskip

{\bf 1.7.A. Further obstructions?} 
Do there exist
further restrictions on the homotopy classes of smooth strictly ergodic
loops in the group of Hamiltonian diffeomorphisms?
I do not know the answer
even in the simplest case when $Y$ is the $2$-sphere
endowed with an area form. In this case the group of Hamiltonian
diffeomorphisms
has the homotopy type of $SO(3)$, and
thus its fundamental group equals $\Z_2$. It would be interesting
to understand
whether in this situation
there exists a smooth strictly ergodic loop in the non-trivial homotopy class.
Note that the obstruction provided by Theorem 1.3.A cannot be applied
since the homotopy class in question has finite order.

\medskip

{\bf 1.7.B. Continuous vs. smooth.} 
Question 1.1.A still makes sense if one considers continuous
loops of Hamiltonian diffeomorphisms instead of smooth ones.
In this situation the existence result 1.2.A above can be refined as
follows. One can show existence of contractible continuous strictly ergodic
loops {\it with every given} irrational rotation number $\alpha$. Note that
in the smooth case the methods used below 
lead to those $\alpha$'s only
which admit a very fast approximation by rationals.
\footnote
{
This feature of the classical approach (see \S 2 below)
is well known to experts. 
It sounds likely however that using methods
developed by Eliasson [E] one can construct
strictly ergodic Hamiltonian skew products on $S^1 \times S^2$
whose rotation
numbers satisfy a diophantine condition.}
On the other hand,
our proof of the obstruction 1.3.A above
does not go through when a strictly ergodic loop is continuous,
since this crucially uses existence of the Hamiltonian function.

\medskip

{\bf 1.7.C. The volume-preserving case.}
Let $Y$ be a closed manifold endowed with a volume form
and let $G$ be the identity component of the group of all volume-preserving
diffeomorphisms. Exactly as in the Hamiltonian case, one can address 
the question about homotopy classes represented by strictly ergodic loops.
The formulation and the proof of the existence result 1.2.A 
remain valid without any changes in the volume-preserving category.
However strictly ergodic loops may well be non-contractible.
We present here an example
of a non-contractible strictly ergodic loop in the group of area
preserving diffeomorphisms of the 2-torus. Note that in dimension 2
every area preserving diffeomorphism is symplectic but not necessarily
Hamiltonian.
This leads to a suggestion that the phenomenon described
in 1.3.A is a purely Hamiltonian one. 

Take 
irrational numbers $\alpha$ and $\beta$ such that $1, \alpha$ and
$\beta$ are independent
over $\Q$.
Consider the loop of transformations of the 2-torus
$h(t):\T^2 \to \T^2$
which take $(y_1,y_2) \in \T^2$ to $(y_1+t,y_2+\beta)$.
Clearly these transformations
preserve the area form $dy_1 \wedge dy_2$
on $\T^2$, and the loop is not contractible.
The standard harmonic analysis argument shows that
the corresponding skew product
$T_{h,\alpha} : \T^3 \to \T^3$ is ergodic. It follows from a theorem due to
Furstenberg [F2, p. 66] that in this case $T_{h,\alpha}$ is strictly ergodic.

The contrast between Hamiltonian and volume-preserving cases
becomes even more transparent when one considers the group 
of volume-preserving transformations
of $Y=S^1$. 
Namely
Furstenberg proved in [F1]
that in this situation {\it every} non-contractible loop
$h : S^1 \to S^1 \subset {\rm Diff}(Y)$ is strictly ergodic.

\bigskip

{\bf 2. Constructing strictly ergodic skew products}

\bigskip

Let $\mu$ be the canonical measure on $Y$.
We write $\cal H$ for the space of all continuous
functions on $Y$ with the zero mean with respect to $\mu$. This space
is endowed with a norm
$||H|| = \max_{y\in Y}|H(y)|$. 
Recall that our task is to prove theorem 1.2.B
on the existence of smooth strictly ergodic loops.

\medskip

{\it 2.1 The classical approach}

\medskip

The proof of 1.2.B is based on the following chain of statements.

\medskip

{\bf Property 2.1.A.} For every continuous function
$F:S^1 \times Y \to \R$ with the zero mean with respect to $dtd\mu$,
for every $\eps > 0$ and for every rational number $r$ there exists a 
loop $g \in C^{\infty}(S^1,G)$
such that the following two conditions hold:

\NI
(i) $|\int_0^1F (t,g(t)^{-1}y)dt| < \eps$ for all $y \in Y$;

\NI
(ii) $g(t+r) = g(t)$ for all $t$. 

\medskip

{\bf Averaging property 2.1.B.} For every $H \in \cal H$ and
$\eps >0$ there exist transformations $g_1,...,g_N \in G$ such that
$${1\over N}|H(g_1^{-1}y) +...+H(g_N^{-1}y)| < \eps,$$
for all $y \in Y$.

Note that 2.1.B is a natural discrete version of 2.1.A(i).
However  in contrast with 2.1.A we
consider here functions $H$ of the variable $y$ only,
and do not care about 
the commutativity condition 2.1.A(ii).

Property 2.1.A implies the statement of Theorem 1.2.B
(see [FH]). Averaging property 2.1.B implies property
2.1.A (see [N] where an analogous implication is proved
in the context of ergodicity; in our situation the argument goes
through without any
essential modifications). For the reader's convenience, we present details
in the Appendix in \S 5.

At this point we face a difficulty. The analogue of 2.1.B for the
$L_1$-norm on $\cal H$, which is used in [N], was proved earlier
by M. Herman [He] with a very elegant use of the Hahn-Banach theorem.
I was unable to adjust Herman's short argument to the $L_{\infty}$-case, 
and thus was forced to take a different route. The key idea is to
derive the averaging property from a certain covering property
which we are going to describe now.

Let us introduce the following useful
object.
Denote by $\cal S$ the set of linear operators $\cal H \to \cal H$
which consists of all averaging operators of the form
$$S^{g_1,...,g_N} (H) =
{1\over N}(H \circ g_1^{-1} +...+H \circ g_N^{-1}),$$
where $N \in \N$ and $g_1,...,g_N \in G$. Note that $\cal S$ is
closed under composition of operators. With this notation 
2.1.B states that for all $H \in \cal H$ there exists $S \in \cal S$
such that $||S(H)||$ is arbitrarily small. An important (and obvious)
feature of transformations from $ \cal S$ is that they do not increase
the norm of functions: $||S(H)|| \leq ||H||$.

\medskip

{\it 2.2 A covering property}

\medskip

{\bf Covering property 2.2.A.} There exist constants $c_1 \geq 0,c_2 \geq 1$
such that for every non-empty
open subset $A \subset Y$ one can find transformations
$g_1,...,g_N \in G$ so that the sets $g_1(A),...,g_N(A)$ form a covering of $Y$
which satisfies the following inequality:
$$ {1\over N}\Sigma_{i=1}^N \chi^{g_i(A)}(y) \geq (c_1 + c_2 {{\mu(Y)}\over
{\mu(A)}})^{-1},$$
for all $y \in Y$.

\medskip

Here and below $\chi^B$ stands for the characteristic function of
a subset $B$.

\medskip

\proclaim Theorem 2.2.B. The covering property implies the 
averaging property.

Note that averaging property 2.1.B 
applied to the normalized characteristic function of an open subset $A$
implies up to $\eps$ covering property 
2.2.A with the optimal constants $c_2 = 1$ and $c_1=0$.
Hence a surprising feature of Theorem 2.2.B is that
starting from an arbitrary choice of the constants
we get the optimal constants. 
Let us mention also that our covering property is motivated
by the Glasner-Weiss covering property [GW].

The rest of the paper is organized as follows. 
Theorem 2.2.B admits a rather short proof which we
present in  section 3. Thus it remains to verify
that
the group of 
Hamiltonian diffeomorphisms enjoys covering property
2.2.A. Here is the idea of our proof. ``Represent'' the symplectic
manifold $Y$ as a cubical polyhedron consisting of small
symplectically standard pairwise equal closed cubes. 
There exists a universal constant, say $k$ which depends
only on $Y$ (but not on the size of the cubes!)
such that $Y$ can be decomposed as the union
of subpolyhedra $Y_1,...,Y_k$ where each $Y_i$ consists
of cubes with {\it pairwise disjoint closure}. Assume
without loss of genericity that the set $A$ given in 2.2.A
is a subpolyhedron of $Y$, and $\mu(A \cap Y_1) \geq \mu(A)/k$.
Set $A_1 = A \cap Y_1$, and assume that this set consists
of $m$ cubes. 
Assume for simplicity that each $Y_i$ consists of $M$ cubes
with $M > m$.
Clearly, every subpolyhedron of $Y_i$ which consists of $m$
cubes is Hamiltonian diffeomorphic to $A_1$. 
Denote by $r$ the number of all such subpolyhedra in $Y_i$,
thus their total number is $N = kr$.
Note that every point of $Y_i$ belongs to at least
$rm/M$ subpolyhedra from our collection.
Thus we can choose $N$ elements of $G$ such that
for every point of $Y$ the left hand side of the inequality  2.2.A
is at least $$rm/Mkr = m/Mk \geq \mu(A)/k\mu(Y).$$
Since $k$ is a universal constant, we get 2.2.A.
The details of this elementary argument are quite cumbersome. They
are worked out  
in section 4 with the use of Katok's results [K].
This will complete the proof of Theorem 1.2.B.

\bigskip

{\bf 3. Covering implies averaging}

\bigskip

In this section we prove theorem  2.2.B.

\medskip

{\it 3.1 A recursive procedure}

Suppose that covering property 2.2.A holds. 
We have to show that
for all $H \in \cal H$ there exists $S \in \cal S$
such that $||S(H)||$ is arbitrarily small.

We claim that it suffices to show that there exists $S_+ \in \cal S$
such that $\max S_+(H)$ is arbitrarily small. Indeed, note that
the same result applied to the function $-S_+(H)$ would show that
there exists $S_- \in \cal S$ such that $\min S_-(S_+(H))$ is arbitrarily
small. But obviously operators from $\cal S$ {\it do not increase}
maximal values of functions, thus $\max S_-(S_+(H)) \leq \max S_+(H)$.
So taking $S = S_- \circ S_+$ we get that $||S(H)||$ is arbitrarily small.
The claim follows.

Assume without loss of generality that $||H|| \leq 1$ and $\mu(Y) = 1$.
We construct the operator $S_+$ with the help of the following
recursive procedure. We start with the function 
$H^{(0)} = H$, and define $H^{(i+1)}$ as the image of $H^{(i)}$
under some specially chosen operator $S_i \in \cal S$. Namely consider
a subset $$A = \{H^{(i)} < {1\over 2} \max H^{(i)}\} \subset Y.$$
Take $g_1,...,g_N$ from the definition of the covering property 2.2.A applied
to the set $A$, and set $S_i = S^{g_1,...,g_N}.$

\proclaim Lemma 3.1.A. The following inequality holds:
$$\max H^{(i+1)} \leq \max H^{(i)} (1- {{\max H^{(i)}\over c}} ),$$
where $c = 2(3c_2 + c_1)$.

The theorem easily follows from the lemma. Set $m_i = \max H^{(i)}$.
Notice that $m_0 \leq 1$ due to our assumption, and the sequence $\{m_i\}$
is non-increasing since operators from $\cal S$ do not increase maximal values
of functions. 
Clearly the sequence $\{m_i\}$ converges to a non-negative number $m$.
For the
proof of the theorem it suffices to show that $m=0$. Assume on the contrary
that $m >0$. Lemma 3.1.A implies that $m \leq m(1-{m\over c})$, which is 
obviously impossible. This contradiction proves the theorem.
\QED

\medskip

{\it 3.2  Proof of 3.1.A} 

\medskip

For simplicity of notations we write $H$ for $H^{(i)}$
and $H'$ for $H^{(i+1)}$. Set $m = \max H$ and $m' = \max H'$.
The proof is divided into 2 parts.

1) Our first task is to find a lower bound for $\mu(A)$.
Since $H$ has zero mean we have
$$\int_A Hd\mu + \int_{Y-A}Hd\mu = 0.$$
The first summand is not less than $-\mu(A)$ since $||H|| \leq 1$
due to our assumption. 
The second summand is not less than ${1\over 2}m(1-\mu(A))$ in view of the
definition of $A$ and our convention that $\mu(Y) = 1$.
Thus $0 \geq -\mu(A) + (1-\mu(A))m/2$ and hence
$$\mu(A) \geq m/(m+2).$$

2) Return now to the definition of $H'$.
Take a point $y \in Y$ and denote by $N'$ the cardinality of the
set $\{j|\; y \in g_j(A)\}$. Clearly,
$$H'(y) \leq {1\over N} ((N-N')m + N'm/2) = m(1-{N'\over {2N}}) .$$
Since this holds for all $y$, the maximum $m'$ of $H'$ satisfies
the same inequality. On the other hand the covering property 2.2.A 
implies that
$$N'/N \geq (c_1 +c_2/\mu(A))^{-1}.$$
Substituting the inequality for $\mu(A)$ obtained in part 1, we get that
$$ N'/N \geq m/(2c_2 + (c_1+c_2)m) \geq m/c_3,$$
where $c_3 = 3c_2 +c_1$. In the last inequality we used that $m \leq 1$.
Finally, substituting this estimate for $N/N'$ into the estimate for $m'$
obtained above, we get that
$m' \leq m(1-m/c)$. This completes the proof.
\QED

\bigskip

{\bf 4. Proving the covering property}

\bigskip

{\it 4.1 Statement of the result}

\medskip

In this section we prove the following result. As before,
let $G$ be the group of Hamiltonian diffeomorphisms of
a closed symplectic manifold $(Y,\Omega)$.

\proclaim Theorem 4.1.A. The group $G$ enjoys covering property 2.2.A.

\medskip

As it was explained in \S 2, this result completes the proof of Theorem 1.2.B.

\medskip

{\it 4.2 A local version}

\medskip

First of all we prove an analogue of
4.1.A for domains
of the linear space $\R^n$. We assume that $\R^n$ is endowed 
with the standard symplectic form.
Let $X \subset \R^n$ be a closed bounded connected domain with piece-wise smooth
boundary.
Let $U$ be an
open domain with compact closure which contains $X$, and write $G_U$ for
the group of all Hamiltonian
diffeomorphisms generated by Hamiltonian functions supported in $U$.
We denote by $\mu$ the canonical measure on $\R^n$. Also, given a family
$\theta = \{A_1,...,A_N\}$ of subsets of $U$, we write $\nu_{\theta}$
for the "counting function" $\Sigma_{i=1}^N \chi^{A_i}$.

\proclaim Proposition 4.2.A. There exist universal constants
$C_1 \geq 0, C_2 \geq 1$ such that for every non-empty open subset $A$ of the
interior of $X$ there exist transformations $g_1,...,g_M \in G_U$
such that $\theta = \{g_1(A),...,g_M(A)\}$ is a covering of $X$, and
$${1\over M}\nu_{\theta}(x) \geq (C_1 + C_2\mu(X)/\mu(A))^{-1},$$
for all $x \in X$.
\medskip

The proof is based on two auxiliary statements.

\proclaim Lemma 4.2.B. For every $a \in (0;\mu(X)]$ there exists a family
$\sigma = \{A_1,...,A_N\}$ of open subsets of $U$ with the compact closure
in $U$ such that $\mu(A_j) = a$ and $\nu_{\sigma}(x)/N \geq a/2\mu(X),$
for all $x \in X$.

\medskip

\proclaim Lemma 4.2.C. There exists a constant $C \in \N$ which depends
only on the dimension of $X$
with the following property.
Let $A,B \subset U$ be two open subsets
with compact closure in $U$ such that $\mu (A) > 2C\mu(B)$.
Then there exist $C$ transformations $g_1,...,g_C \in G_U$ such that
$ B \subset \cup_{i=1}^C g_i(A)$. 

Let us derive Proposition 4.2.A from these lemmas.

{\it Proof of 4.2.A:}
Take any open subset $A$ as in 4.2.A and apply Lemma 4.2.B with
$a = \mu(A)$. We get a family $\sigma = \{A_1,...,A_N\}$ of subsets.
Fix $\eps >0$. We claim that there exist $f_1,...,f_N \in G_U$
such that $\mu (f_i(A) \Delta A_i) < \eps/N$ for all $i=1,...,N$.
The claim follows immediately from Katok's Basic Lemma [K].
Consider a family of subsets
$\tau = \{f_1(A),...,f_N(A)\}$, and set $B = \cup_{i=1}^N (f_i(A) \Delta A_i).$
Clearly, $\mu(B) < \eps$.
Assume now that $\eps$ is so small that Lemma 4.2.C can be applied 
to $A$ and a small open neighbourhood of $B$. 
Using this lemma, we get $C$ transformations
$g_1,...,g_C \in G_U$ such that $ B \subset \cup g_i(A)$. 

Set $N'= \min_{x \in X} \nu_{\sigma}(x)$. Define a new family
$\theta$ which consists of all subsets of the form $g_i(A), i= 1,...,C$
taken $N'$ times, and in addition of all subsets from $\tau$. In other words,
$$\theta = \{f_1(A),...,f_N(A);g_1(A),...,g_C(A),...,g_1(A),...,g_C(A)\}.$$
The number $M$ of elements in $\theta$ equals to 
$N + CN'$. On the other hand, we claim that $\nu_{\theta}(x) \geq N'$
for all $x \in X$. Indeed, for $x \in B$ this follows from the
definition of transformations $g_i$. For $x \in X-B$ we note that
$\nu_{\tau}(x) = \nu_{\sigma}(x)$ and the claim follows from the
definition of $N'$. 

Recall now from Lemma 4.2.B that $N'/N \geq \mu(A)/2\mu(X)$.
Thus
$${1\over M}\nu_{\theta} \geq N'/(N+CN') \geq (C+2\mu(X)/\mu(A))^{-1}.$$
Therefore we proved 4.2.A with $C_1 = C$ and $C_2 = 2$.
\QED

\medskip

It remains to prove the lemmas. Both are of combinatorial nature,
and we need some suitable notions. By {\it a cubical partition of size} $u$
we mean a decomposition of $\R^n$ into equal closed cubes of volume $u$.
The cubes may intersect along the boundaries only, and their centers form
a lattice. 
{\it A cubical polyhedron} is the union 
of some cubes of a cubical partition.

\medskip

{\it Proof of 4.2.B:} 
Take a sufficiently large positive integer $k$ and consider a cubical
partition of size $a/k$. We can assume that there exists 
a cubical polyhedron $X', \;\; X \subset X' \subset U$ such that
$\mu(X') < 1.5\mu(X)$ and $\mu(U-X') > a/k$.
Take all possible
sub-polyhedra  of $X'$ consisting of exactly $k-1$ cubes. Let $A_1,...,A_N$
be their open neighbourhoods of volume $a$ which have compact closure in $U$.
Denote this family by $\sigma$.
Clearly for all $x \in X$
$$\nu_{\sigma}(x)/N \geq (k-1)a/k\mu(X') \geq a/2\mu(X).$$
This proves the lemma.
\QED

\medskip

For the proof of the second lemma we need the following facts
from elementary geometry of $\R^n$. 
Given a closed cube $Q \subset \R^n$
and a positive number $c$ denote by $cQ$ the cube with the same center
which is homothetic to $Q$ with the coefficient $c$. 
Let $Q_1,Q_2$ be two cubes which are obtained from each other
by a translation. Note first that
if $Q_1\cap Q_2 \neq \emptyset$ then $Q_1 \subset {\rm Interior}(4Q_2).$
Define
a {\it nice subpartition} of a cubical partition as a family of pairwise disjoint cubes
from the partition  which satisfies the following property. Given any two cubes
$Q_1$ and $Q_2$ from the family, their homothetic images $16Q_1$ and
$16Q_2$ are also disjoint.
Note finally that there exists a constant $C \in \N$ which depends only on $n$
such that every cubical partition of $\R^n$ can be decomposed into $C$
nice subpartitions.

\medskip

{\it Proof of 4.2.C:}

1) Let $U' \subset U$ be an open connected
domain with compact closure in $U$ which
contains both $A$ and $B$. 
Fix a cubical partition $P$ of size $u$.
We assume that $u$ is so small that for each cube $Q$ of the partition
which intersects $U'$ holds $16Q \subset U$.
Let $P = P_1 \cup ...\cup P_C$ be its decomposition into nice subpartitions.
Denote by $A_i$ the union of all cubes from $P_i$ which are contained in $A$,
and by $B_i$ the union of all cubes from $P_i$ which intersect $B$.
Clearly taking $u$ small enough we can achieve that
$\mu(B_i) \leq 1.1\mu(B)$ for all $i$, and $\mu(A_i) \geq \mu(A)/1.1C$
for some $i$. Assume without loss of generality that the last inequality
holds for $i=1$. 
Since $\mu(A)> 2C\mu(B)$ we conclude that the number of cubes in $A_1$
is greater than the number of cubes in each of $B_i$. 
Fix some $i \in \{1,...,C\}.$
Clearly, in order to prove the lemma it suffices to show that
if $u$ is small enough there exists a transformation
$g_i \in G$ such that $B_i \subset g_i(A_1)$. 

\medskip

2) Denote by $q =\{Q_1,...,Q_k\}$ the set of all cubes in $A_1$,
by $q'=\{Q_1',...,Q_m'\}$ the set of all cubes in $B_i$,
and by $p$ the union $q \cup q'$. We write $4q$ for 
$\{4Q_1,...,4Q_k\}$, and define analogously $4q'$ and $4p$.
Let $Z$ be the subset of $U$ obtained by the union of all cubes
from $4p$. We claim that each connected component of $Z$ is either
one cube, or the union of two cubes. Assume on the contrary that
there exist three cubes from $4p$ such that two of them intersect
the third one. By definition of a nice subpartion, either
these two belong to $4q$ and the third belongs to $4q'$, or vice
versa the two belong to $4q'$ and the third to $4q$. Without loss 
of genericity we assume that $4Q_1$ and $4Q_2$ intersect $4Q_1'$.
But then $16Q_1$ and $16Q_2$ contain $4Q_1'$, and that contradicts
the definition of a nice subpartition. The claim follows.

\medskip

3) In view of step 2, we may assume without loss of generality that
for some $r \leq m$ 
$$4Q_1 \cap 4Q_1' \neq \emptyset,..., 4Q_r \cap 4Q_r' \neq \emptyset,$$
and all other pairs of cubes from $4p$ are disjoint. Recall also that
$m < k$. We claim that for every $j \in \{1,...,m\}$ there exists
a transformation $h_j \in G_U$ such that $h_j(Q_j)  = Q_j'$ and
$h_j$ equals the identity on all $Q_l,Q_l'$ with $l \neq j$.
Note that this claim implies the lemma. Indeed set
$g_i = h_1 \circ...\circ h_r$. Clearly $g_i(Q_j) = Q_j'$,
and hence we constructed a transformation as required in step 1.

\medskip

4) It remains to prove the claim of step 3. Take $j \in \{1,...,m\}$.

First assume that $j \leq r$. Consider a set $K = 16Q_j \cap 16Q_j'$.
Clearly, $K$ is a convex polyhedron
whose interior contains both $Q_j$ and $Q_j'$.
There exists a path of transformations from $G_U$
supported in $K$ whose time one map takes $Q_j$ to $Q_j'$. 
Take this time one map as $h_j$.
Any other cube from $p$ is disjoint from $K$ by definition of a nice
subpartition. 
Thus $h_j$ has the required properties.

Assume now that $j > r$. Consider the set $Z_j = Z - (4Q_j \cup 4Q_j')$.
It follows from step 2 that $U' - Z_j$ is a connected set. Join the
centers of cubes $Q_j$ and $Q_j'$ by a smooth path $\gamma \subset U' - Z_j$.
Let $t$ be a parameter along $\gamma$ which runs from $0$ to $1$.
Denote by $K_t$ the cube centered in $\gamma(t)$ which is obtained from
$Q_j$ by a parallel translation.
Since the center of $K_t$ is disjoint from all cubes $4Q_l,4Q_l'$ with 
$l \neq j$
then $K_t$ is disjoint from all such $Q_l$ and $Q_l'$. In particular,
there exists a small neighbourhood say $V$ of the union of all cubes
$K_t$ which is disjoint from all cubes $Q_l,Q_l'$ with $l \neq j$.
It is easy to see that there exists a path of transformations
from $G_U$ supported in $V$ whose time one map takes $Q_j$ to $Q_j'$.
Take this time one map as $h_j$. This completes the proof.
\QED

\medskip

{\it 4.3 Proof of 4.1.A}

\medskip

Denote by $\mu$ the canonical measure on $Y$.
Consider a triangulation $X_1,...,X_r$ of $Y$ such that every simplex
$X_i$ is contained in a (Darboux) coordinate chart $U_i$. Assume moreover
that all $X_i$ have equal volume $\mu(Y)/r$.
Let $A \subset Y$ be an open subset. A straightforward application
of Katok's Basic Lemma [K] 
shows that there exists
a transformation $f \in G$ with the following property: 
$$\mu(f(A)\cap X_i) > \mu(A)/2r,$$
for all $i = 1,...,r$.

Obviously it suffices to check the covering property for 
the set $ A' = f(A)$. Let $A_i$ be an open set which lies in the interior
of $A' \cap X_i$ and has volume $\mu(A)/2r$.
Apply the local statement 4.2.A to the triple $(U_i,X_i,A_i)$.
For every $i$, we get a sequence of transformations $g_{ij},\; j=1,...,N_i$
such that for all $x \in X_i$ holds

$${1\over N_i} \Sigma_{j=1}^{N_i} \chi^{g_{ij}(A_i)} (x) \geq
(C_1 + C_2\mu(X_i)/\mu(A_i))^{-1} =  \lambda, $$
where
$$\lambda = (C_1 +2C_2\mu(Y)/\mu(A'))^{-1}.$$
Repeating the sequences we can achieve that each of them has the same number
of terms, in other words that all $N_i$'s are equal to the same number $N$.
We claim that the family of transformations
$$\{g_{ij}\}, \; i= 1,...,r \; ; \; j = 1,...,N$$
does the job for the set $A'$ with universal constants
$c_1 = rC_1$ and $c_2 = 2rC_2$.
Indeed, take $y \in Y$. Without loss of generality assume that
$y \in X_1$. Thus
$${1\over {Nr}} \Sigma_{i,j}\chi^{g_{ij}(A')} (y) \geq
{1\over {Nr}}\Sigma_{j} \chi^{g_{1j}(A_1)} (y) \geq {{\lambda N} \over {rN}}=
(rC_1 +2rC_2\mu(Y)/\mu(A'))^{-1}.$$
This completes the proof.
\QED


\vfill\eject

{\bf 5. APPENDIX: more details on the classical approach}

\bigskip

In this appendix we present details of the classical approach
to constructing strictly ergodic Hamiltonian skew products
(see section 2.1 of the main text). 

\medskip

{\it 5.1 Averaging property 2.1.B implies 2.1.A} (cf. [N])

\medskip

We work in the notations of section 2.1. The proof is divided
into several steps.


1) We claim that for every finite sequence $H_1,...,H_k$ of functions
from $\cal H$ and for every $\eps >0$  there exists an operator
$S \in {\cal S}$ such that $||S(H_i)|| < \eps$ for all $i$.
Indeed, using averaging property 2.1.B define recursively a
sequence of operators $S_i \in {\cal S}$ such that the following inequalities
hold:
$$||S_1(H_1)|| < \eps, \; ||S_2(S_1(H_2))|| < \eps, ... , 
||S_k(...(S_1(H_k)...)||
< \eps.$$
Set $S = S_k \circ ... \circ S_1$. The operator $S$ is as required in view
of the fact that operators from $\cal S$ do not increase the norm of functions.

As an immediate consequence of the claim we get that for every sequence
$H_1,...,H_k$ and every $\eps >0$ there exists a smooth loop $h:S^1 \to G$ such that
$$|\int_0^1 H_i(h(t)^{-1}y) dt | < \eps,$$
for all $y \in Y$.

\medskip

2) Note that in order to verify condition 2.1.A(i) we can assume
that 
for each fixed $t$
the function $F(t,.)$ belongs to the space $\cal H$. Indeed,
given any $F$ 
one can modify
it as follows:
$$F'(t,y) = F(t,y)- \int_Y F(t,z) d\mu(z).$$
Clearly our assumption holds for $F'$.  Moreover $F$ and $F'$ satisfy or do not
satisfy the conditions in question simultaneously. 

\medskip

3) Let $F$ be a function which satisfies the assumption
of step 2.
We claim that for every $\eps > 0$
there exists a smooth loop $h$ such that for all $s \in S^1$
the following inequality holds:
$$|\int_0^1 F(s,h(t)^{-1}y) dt | < \eps/3.$$
Indeed, choose a large natural number $N$ such that
$|F(t',y)-F(t'',y)| < \eps/9$ for all $y$ provided $|t'-t''| < 1/N$.
Set $p_i = i/N$ where $i=0,...,N-1$. Using the last statement of step 1 choose
a smooth loop $h(t)$ such that 
$$|\int_0^1 F(p_i, h(t)^{-1}y) dt | < \eps/9.$$
A straightforward estimate shows that $h$ is as required.

\medskip

4) Take now an arbitrary integer $M > N$. A more specific choice of $M$
will be made in the next step. Meanwhile we denote $g(t) = h(Mt)$
and claim that 
$$ I = |\int_0^1 F(t,g(t)^{-1}y) dt | < \eps.$$
Here is the proof. Set $q_i = i/M$ where $i=0,...,M$. An obvious argument
shows that
$$I \leq \eps/3 + \Sigma_{i=0}^{M-1}|J_i|,$$
where $$J_i= \int_{q_i}^{q_{i+1}} F(q_i,h(Mt)^{-1}y)dt.$$
Introducing a new variable $s = Mt-i$ we get that
$$J_i = {1\over M} \int_0^1 F(q_i,h(s)^{-1}y)ds,$$
and thus our choice of $h$ implies that $|J_i| < \eps/(3M)$ for all $i$.
The claim follows immediately.

\medskip

5) Let us sum up the results of the previous steps. We constructed 
a loop $g(t) = h(Mt)$ which satisfies 2.1.A(i) above.
Moreover the choice of the sufficiently large integer $M$ is in our hands.
We are going to use this in order to guarantee condition 2.1.A(ii).
Let $r$ be a rational number.
Taking $M$ as a large multiple of the denominator of $r$ we get
that $g(t+r)=g(t)$. 

This completes the proof.
\QED


\medskip

{\it 5.2 Property 2.1.A implies Theorem 1.2.B}

\medskip

The proof of this statement occupies the rest of the appendix.
We follow very closely the exposition in [FH] (see also [GW]).
Let us fix some notations.
We work on the manifold $X =S^1 \times Y$. Denote by $D$ the group
of all skew products $(t,y) \to (t+\alpha,h(t)y)$ where $\alpha \in S^1$
and $h :S^1 \to G$ is a smooth loop. Let $D_0$ be its subgroup consisting
of maps of the form $(t,y) \to (t,h(t)y)$. We write $S_{\alpha}$ for the
shift $(t,y) \to (t+\alpha,y)$. 

Consider the set of all mappings of the form 
$\phi^{-1}\circ S_{\alpha} \circ \phi$
where $\phi \in D_0$ and $\alpha \in S^1$. Its closure in $D$ can be identified
in an obvious way with the set $\cal R$ introduced in 1.2, and in this appendix
we use for it the same notation $\cal R$.

\medskip

{\it 5.3 Minimality and unique ergodicity}

\medskip

Consider two collections $\cal A$ and $\cal B$ of subsets of $\cal R$
as follows. The collection $\cal A$ consists of all subsets of the
form
$$A(U) = \{T \in {\cal R}| \cup_{i =0}^{\infty} T^iU = X\},$$
where $U$ runs over all non-empty open subsets of $X$.
The collection $\cal B$ consists of all subsets of the form 
$$B(F,\eps) = \{T \in {\cal R}| \inf_{N\in \N} ||{1\over N}
\Sigma_{j=0}^{N-1} F \circ T^j || < \eps\},$$
where $F$ runs over all continuous functions with zero mean on $X$,
and $\eps > 0$.
Note that the intersection of all sets from $\cal A$ consists of minimal
diffeomorphisms (that is every orbit is dense in $X$). Also, the intersection
of all sets from $\cal B$ consists of uniquely ergodic diffeomorphisms
(that is there exists precisely one invariant Borel probability measure).
Write $\cal C$ for the
union of $\cal A$ and $\cal B$. 
It is easy to see that minimality combined with unique ergodicity
implies strict ergodicity. Thus for existence of strictly
ergodic skew products it suffices to show that the intersection of
all sets from $\cal C$ is non-empty. 
Consider {\it a countable}
subcollection ${\cal C}'$ of
$\cal C$
which consists of all sets of the form $A(U_i)$ and $B(F_j,{1\over k})$,
where $\{U_i\}$ is a countable basis of open subsets on $X$, $\{F_j\}$
is a countable dense subset of the space of continuous functions
with zero mean on $X$, and the number $k$ runs over the natural numbers.
Obviously every set from $\cal C$ contains a set from ${\cal C}'$,
thus the intersection of all sets from ${\cal C}$ is equal to the
intersection of all sets from ${\cal C}'$.

Now the strategy is as follows. The collection $\cal C$ consists
of subsets which are open in $\cal R$. Since the group $G$ is closed
in ${\rm Diff}(Y)$,
the set $\cal R$ (with the topology induced from ${\rm Diff}(X)$)
has the Baire property.
\footnote
{This is exactly the place where we use that $G$ is closed
in ${\rm Diff}(Y)$. It seems however that one can 
prove the theorem without this assumption. For that purpose
one should work with a different cleverly chosen topology on $G$.}
 Therefore it suffices to show that
the subsets from $\cal C$ are dense in $\cal R$.

\proclaim Lemma 5.3.A. Let $C$ be a set from ${\cal C}$ and
$r$ be a rational
number. 
There exists $\phi \in D_0$ which commutes with $S_r$ and
such that for every irrational number $\alpha$ the
diffeomorphism $\phi^{-1} \circ S_{\alpha}\circ \phi$ belongs to $C$.

\medskip

{\it 5.4 The final argument}

\medskip

Assume the Lemma. We claim that for every rational $r$ and every
$C\in {\cal C}$ the shift $S_r$ belongs to the closure of $C$.
Indeed, take $\phi$ from the Lemma and
choose a sequence $\{\alpha_j\}$ of irrational numbers
which converges to $r$. 
Then 
$\{\phi^{-1}\circ S_{\alpha_j}\circ \phi\}$ is a sequence of elements
of $C$ which converges to $S_r$, and the claim follows.
Since rational numbers are  dense in the circle, we conclude that every
(rational or irrational) shift is in the closure of $C$.
Note that for
every $\psi \in D_0$ and $C \in {\cal C}$ the set 
$$ \{\psi \circ f \circ \psi^{-1} | f \in C\}$$
is again contained in $\cal C$. Thus every element of 
the form $\psi^{-1} \circ S_{\alpha} \circ {\psi}$ with
$\psi \in D_0$ belongs to the closure of $C$. Since by definition
these elements are dense in $\cal R$, we get that $C$ is dense
in $\cal R$. As it was explained in 5.3 this completes the proof 
of existence of smooth strictly ergodic skew products.

\medskip

{\it 5.5 Proof of 5.3.A} 

\medskip

Fix a set $C \in \cal C$ and a rational number $r$.
We have to consider two cases.

1) Let $C = A(U)$. We can assume that $U$ splits as $\Delta \times V$, where
$\Delta $ is an open interval of $S^1$ and $V$ is an open subset of $Y$.
Take a loop $g:S^1 \to G$ in such a way that the sets
of the family $\{g_t(V)\},\; t \in \Delta$ cover $Y$. Moreover
it is easy to achieve that $g(t+r)=g(t)$ for all $t$.
Define an element $\phi \in D_0$ by $\phi(t,y) = (t,g(t)y)$. It
commutes with $S_r$. Take an irrational $\alpha$.
Write the union 
$$\cup_{i=0}^{\infty} (\phi^{-1} \circ S_{\alpha}\circ \phi)^i (U)$$ 
as
$\phi^{-1}(W)$ where
$$W= \cup_{i=0}^{\infty}S_{\alpha}^i (\phi(U)).$$
In view of our construction the set $\phi(U)$ intersects the circles
$S^1 \times \{y\}$ for all $y \in Y$. Since 
orbits of $S_{\alpha}$ are dense on each such circle
we conclude that $W=X$. Then $\phi^{-1}(W)=X$ and therefore 
$\phi$ is as needed.

\medskip

2) Suppose now that $C=B(F,\eps)$. Given a loop $g:S^1 \to G$ consider
the integral
$$I(y) = \int_0^1 F(t,g(t)^{-1}y) dt.$$
Using Property 2.1.A we can choose $g$ in such a way that $g(t+r)=g(t)$
for all $t$, and $|I(y)|<\eps/2$ for all $y$.
Define an element $\phi \in D_0$ by $\phi(t,y) = (t,g(t)y)$.  It commutes
with $S_r$. Take an irrational $\alpha$. Write the ergodic sum 
$${1\over N}\Sigma_{j=0}^{N-1} F \circ (\phi^{-1}\circ S_{\alpha}\circ\phi)^j $$
as $ G_N \circ \phi$ where
$$G_N = 
{1\over N}\Sigma_{j=0}^{N-1} (F \circ\phi^{-1}) \circ S_{\alpha}^j .$$
Since the shift of the circle  $t \to t+\alpha$ is strictly ergodic,
and the family $\{F(t,g(t)^{-1}y)\},\; y \in Y$
of functions $S^1 \to \R$ is compact, the ergodic sum above converges
uniformly to $I(y)$ when $N$ goes to infinity.
In particular for large $N$ holds $||G_N \circ \phi|| = ||G_N|| < \eps$,
and we conclude that $\phi$ is as needed.
This completes the proof.
\QED



\vfill\eject

{\bf Acknowledgements.} I thank Viktor L. Ginzburg,
Eli Glasner, Dusa McDuff, Vitali Milman,
J\"urgen
Moser, and Paul Seidel
for very useful discussions, consultations and suggestions.
I am very grateful to Felix Schlenk for pointing out a number of
inaccuracies in the first version.
This paper was written
during my stay at ETH-Z\"urich.
I thank the Forschungsinstitut f\"ur Mathematik
for the hospitality.

\bigskip

\NI
{\bf References}
\bigskip

\NI
[AK] D.~Anosov and A.~Katok, New examples in smooth Ergodic Theory.
Ergodic diffeomorphisms, {\it Trans. Moscow Math. Soc.} {\bf 23} (1970),
1-35.
\MS

\NI
[AM]
M.~Abreu and D.~McDuff, in preparation.
\MS

\NI
[BB] D.~Berend and V.~Bergelson, Ergodic and mixing sequences
of transformations, {\it Ergodic Th. Dynam. Syst.} {\bf 4}(1984), 353-366.
\MS

\NI
[BP]
M.~Bialy and L.~Polterovich,
Invariant tori and Symplectic Topology, {\it Amer. Math. Soc. Transl.} 
(2)
{\bf 171} (1996), 23-33.
\MS

\NI
[E] 
L.H.~Eliasson, Ergodic skew-systems on $T^d \times SO(3,\R)$,
Preprint ETH Zurich, 1991.
\MS

\NI
[F1]
H.~Furstenberg, Strict ergodicity and transformations of the torus,
{\it Amer. J. Math.} {\bf 83} (1961), 573-601.
\MS

\NI
[F2]
H.~Furstenberg, {\it Recurrence in Ergodic Theory and Combinatorial 
Number Theory}, Princeton University Press, 1981.
\MS

\NI
[FH]
A.~Fathi and M.~Herman, Existence de diffeomorphisms minimaux,
{\it Asterisque} {\bf 49} (1977), 37-59.
\MS

\NI
[GW]
S.~Glasner and B.~Weiss, On the construction of minimal skew products,
{\it Israel J. Math} {\bf 34}(1979), 321-336.
\MS

\NI
[H]
H.~Hofer, On the topological properties of symplectic maps,
{\it Proc. Roy. Soc. Edinburgh Sect A} {\bf 115}(1990), 25-38.
\MS

\NI
[He] 
M.~Herman, Construction de diffeomorphismes ergodiques, Unpublished manuscript.
\MS

\NI
[K]
A.~Katok, Ergodic perturbations of degenerate integrable Hamiltonian systems,
{\it Math. USSR Izvestija} {\bf 7}(1973), 535-571.
\MS

\NI
[Kr]
U.~Krengel, {\it Ergodic Theorems}, Walter de Gruyter, 1985.
\MS

\NI
[LM]
F.~Lalonde and D.~McDuff, Hofer's
$L^{\infty}$-geometry: energy and stability of Hamiltonian flows I,II,
{\it Invent. Math.} {\bf 122}  (1995), 1-69.
\MS

\NI
[LMP] F.~Lalonde, D.~McDuff and L.~Polterovich, On the Flux
Conjectures, to appear in the CRM Proceedings and Lecture Notes, 1998.
\MS

\NI
[MS] D.~McDuff and D.~Salamon, {\it Introduction to Symplectic Topology},
Clarendon Press, Oxford, 1995.
\MS

\NI
[N] M.~Nerurkar, On the construction of smooth ergodic
skew-products, {\it Ergod. Th. \& Dynam. Sys.} {\bf 8} (1988), 311-326.
\MS

\NI
[P1]
L. Polterovich, Hamiltonian loops and Arnold's 
principle, {\it Amer. Math. Soc. Transl.} (2) {\bf 180} (1997), 181-187.
\MS

\NI
[P2]
L. Polterovich,  Precise measurements in Symplectic Topology,
to appear in {\it Proceedings of the 2-nd European Congress, Budapest, 1996.}
\MS

\NI
[Si]
J.-C.~Sikorav, {\it Systemes Hamiltoniens et topologie symplectique},
ETS Editrice, Pisa, 1990.
\MS

\NI
[U]
I.~Ustilovsky, Conjugate points on geodesics of 
Hofer's metric, {\it Diff. Geometry and its Appl.}
{\bf 6} (1996), 327-342.
\MS

\end{document}